\documentclass[preprint,12pt]{elsarticle}
\usepackage{amssymb}
\usepackage{amsmath,amssymb}
\newtheorem{thm}{Theorem}

\newtheorem{lem}[thm]{Lemma}
\newtheorem{prop}[thm]{Proposition}
\newtheorem{cor}[thm]{Corollary}
\newdefinition{rmk}{Remark}
\newproof{pf}{\textbf{Proof}}
\newcommand{\supp}{\operatorname{supp}}
\journal{Applied Mathematics and Computation}

\begin{document}

\begin{frontmatter}



\title{Frames for compressed sensing using coherence}

 \author[Gavruta]{L. G\u avru\c ta}
 \ead {gavruta\_laura@yahoo.com}
\author[Eskandani]{G. Zamani Eskandani}
 \ead{zamani@tabrizu.ac.ir}
 \author[Gavruta]{P. G\u avru\c ta}
 \ead{pgavruta@yahoo.com}
  \address[Gavruta]{"Politehnica" University of Timi\c soara, Department of Mathematics, Pia\c{t}a Victoriei no.2, 300006 Timi\c{s}oara, Romania}
\address[Eskandani]{Faculty of Sciences, Department of Mathematics, University of Tabriz,Tabriz, Iran}

\begin{abstract}
We give some new results on sparse signal recovery in the presence of noise, for weighted spaces. Traditionally, were used dictionaries that have the norm equal to 1, but, for random dictionaries this condition is rarely satisfied. Moreover, we give better estimations then the ones given recently by Cai, Wang and Xu.
\end{abstract}


\end{frontmatter}


\section{Introduction}
Compressed sensing (also known as compressive sensing or compressive sampling) is a relatively new field of research, started by the work of E. Cand\'es, D. Donoho, J. Romberg and T. Tao (see \cite{Candes2}, \cite{Candes3}, \cite{Donoho1}). Using the theory of compressed sensing we can recover signals and images from far fewer samples or measurements than were traditionally thought necessary. In order to do this, compressed sensing is based on two principles: sparsity and
incoherence. Using the concept of sparsity it is possible to exactly recover a signal $c$ knowing that it is
sparse or nearly sparse in the sense that it has only a limited number of nonzero components. On the other hand, incoherence says that unlike the signal of
interest, the sampling/sensing waveforms have an extremely dense representation in a proper basis (\cite{Candes1}).

Important applications of compressed sensing are in signal processing, imagine processing (\cite{Duarte}), coding and information theory (\cite{Candes2}), compressive radar (\cite{Baraniuk}, \cite{Potter}), MRI (Magnetic Resonance Imaging), where using compressed sensing techniques we can obtain benefits of imagine speed, reducing costs (see for example \cite{Lustig1}, \cite{Lustig2}). Successfully, a significant number of compressive sensing recovery algorithms were discovered, for example orthogonal matching pursuit (\cite{Mallat}), basis pursuit (\cite{Chen}), $l_1$-minimization. Current directions
of research are in computational biology, geophysical data analysis, astronomy, communications and much more other areas.

 In this paper we present a new result on sparse signals recovery in the presence of noise, which generalizes and completes a result of T. Cai et. all, presented in the paper \cite{Cai}. Our paper is in connection to frame theory; we use the synthesis operator to establish new results on mutual coherence.

Let $\mathcal{H}=\mathcal{H}_n$ be an real $n-$ dimensional space, with the inner product $\langle\cdot,\cdot\rangle$ and $\mathcal{F}=\{f_1,f_2,\ldots,f_N\}\subset\mathcal{H}_n.$ We associate the following operators:\\

The \textit{analysis operator} $$\Theta:\mathcal{H}_n\rightarrow \mathbb{R}^N,$$ which is given by $$\Theta x:=(\{\langle x,f_1\rangle,\langle x,f_2\rangle,\ldots,\langle x,f_N\rangle\})$$ and the \textit{synthesis operator} $$T:\mathbb{C}^N\rightarrow\mathcal{H},\quad T(c_1,c_2,\ldots,c_N)=\sum_{j=1}^Nc_jf_j.$$
We suppose that $w_i=\|f_i\|\neq 0,\hspace{2mm}i=\overline{1,N}.$\\

For $c=(c_1,c_2,\ldots,c_N)\in\mathbb{C}^N$ and for $0<p<\infty$ we denote
$$\|c\|_{p,w}:=\bigg(\sum_{i=1}^N|c_i|^pw_i^p\bigg)^{1/p},\quad\langle c,d\rangle_w=\sum_{i=1}^N c_id_iw_i^2$$ and
$$\langle c,d\rangle=\sum_{i=1}^N c_id_i.$$ We define by $$\mu:=\max_{\substack{i\neq j}}\frac{|\langle f_i,f_j\rangle|}{\|f_i\|\|f_j\|}$$ the coherence of $\mathcal{F}.$  Usually, it is assumed that $\|f_i\|=1, i=\overline{1,N}$
but, for random dictonaries is very rarely satisfied.
It is known that  $\mu$ satisfies the Welch's inequality \cite{Welch}$$\sqrt{\frac{N-n}{n(N-1)}}\leq\mu\leq 1$$
We also denote by $\|c\|_0=\#\{i:c_i\neq 0\}$ the cardinality of the support of $c$.\\
Clearly, $\|c_1+c_2\|_0\leq\|c_1\|_0+\|c_2\|_0$, but $\|\cdot\|_0$ isn't homogeneous.
 We say that $c$ is $s-$sparse if $\|c\|_0\leq s.$

It is known that if $\mathcal{F}$ is a frame for $\mathcal{H}_n$, then the equation $$y=Tc$$ has a solution for any $y\in\mathcal{H}_n.$

We consider the equation $$y=Tc+z$$ where $z$ is an unknown noise term, with $\|z\|_2\leq\varepsilon$ and the problem $$(P_{1,w})\quad\min_{\tilde{c}\in\mathbb{R}^N}\|\tilde{c}\|_{1,w}\quad subject\hspace{1.5mm}to\hspace{1.5mm}\|y-T\tilde{c}\|_{2,w}\leq\eta.$$
The solution of this problem is given in the final part of the present paper. Our result extends and gives better estimations than the one presented in reference \cite{Cai}. Also, we present a result related to Orthogonal Matching Pursuit algorithm (in connection with a result given by J.A. Tropp \cite{Tropp}) for dictionaries with the norm not necessarily equal to $1$.
 \section{Preliminary results}
\begin{lem}(Basic Lemma)\label{Basic lemma} Let $\mathcal{H}_n$ be an real $n$-dimensional Hilbert space and $T$ the synthesis operator for $\mathcal{F}$. Then, for all $c,d\in\mathbb{R}^N$, we have
\begin{enumerate}[$(i)$]
\item $\supp c \cap\supp d=\emptyset$ implies $|\langle Tc,Td\rangle|\leq\mu\|c\|_{1,w}\|d\|_{1,w}$\\
\item  $(1+\mu)\|c\|_{2,w}^2-\mu\|c\|_{1,w}^2\leq\|Tc\|^2\leq(1-\mu)\|c\|_{2,w}^2+\mu\|c\|_{1,w}^2$\\
\item if $c$ is $s$-sparse, then we have
$$ [1-\mu(s-1)]\|c\|_{2,w}^2\leq\|Tc\|^2\leq[1+\mu(s-1)]\|c\|_{2,w}^2$$
\end{enumerate}

\end{lem}
\begin{pf}For all $c,d\in\mathbb{R}^N$ we have
\begin{align*}
|\langle Tc,Td\rangle-\langle c,d\rangle_w|&=\bigg|\langle\sum_{i=1}^N c_if_i,\sum_{j=1}^N d_jf_j\rangle-\sum_{i=1}^N c_id_iw_i^2\bigg|\\
                                           &=\bigg|\sum_{i\neq j}c_id_j\langle f_i,f_j\rangle\bigg|\\
                                           &\leq\mu\sum_{i\neq j}|c_i||d_j|\|f_i\|\|f_j\|\\
                                           &=\mu(\|c\|_{1,w}\|d\|_{1,w}-\langle c,d\rangle_w)
\end{align*}

The relation in $(i)$ follows  from the above computation because $\supp c \cap \supp d=\emptyset$ implies $\langle c,d\rangle_w=0.$

The relation in $(ii)$ follows for $d=c$, and the relation in $(iii)$ it follows from $(ii)$ using Cauchy-Schwarz inequality.
\end{pf}
The following three Propositions are extensions to weighted spaces of some well-known results.
\begin{prop}
Let $\|c\|_0\leq s$, $\|d\|_0\leq s$ and $Tc=Td$. For $$s<\frac{1}{2}\bigg(\frac{1}{\mu}+1\bigg)$$ it follows that $c=d.$
\end{prop}
\begin{pf}
$0=\|T(c-d)\|^2\geq[1-\mu(2s-1)]\|c-d\|_{2,w}^2$. It follows that $\|c-d\|_{2,w}=0$ which implies that $c=d.$
\end{pf}
\begin{prop}
Let $s<1+\frac{1}{\mu}$. Then $\{f_1,\ldots,f_s\}$ is linear independent.
\end{prop}
\begin{pf}
If $$\sum_{i=1}^s c_i f_i=0,$$ then, by Lemma \ref{Basic lemma}, it follows $$\|c\|_{2,w}^2=0$$ hence $c=0.$
\end{pf}
\begin{prop}\label{prop1sup}
Let $$\delta_s=\sup_{\substack{c\neq 0\\\|c\|_0\leq s}}\frac{|\|Tc\|^2-\|c\|_{2,w}^2|}{\|c\|_{2,w}^2}.$$
Then $\delta_s\leq\mu(s-1).$
\end{prop}
\begin{pf}
Immediately from Lemma \ref{Basic lemma}.
\end{pf}

When the elements of $\mathcal{F}$ are all of norm $1$ and $\mu<\dfrac{1}{2s-1}$, J.A. Tropp \cite{Tropp} shown that Orthogonal Matching Pursuit algorithm will recovery any $s-$sparse signal for measurements $y=Tc.$

 We will show that, in the general case when the elements are not necessarily of norm equal to $1$, we have a more general result with a much more easier proof.

 \begin{thm}\label{thm2sup} Let $c\neq 0$, $c\in\mathbb{R}^N$ and $j_0$ be such that $$\langle Tc,f_{j_0}^{'}\rangle=\min_{1\leq j\leq N}\langle Tc,f_{j}^{'}\rangle,\quad where\quad f_j^{'}=\frac{f_j}{\|f_j\|},\quad 1\leq j\leq N.$$
 If $c$ is $s-$sparse and $\delta_s+\mu s<1$ then $j_0\in\supp c.$
 \end{thm}
 \begin{pf}If $j_0\notin\supp c,$ then we have \begin{equation}\label{suportul}\mu\geq\frac{\|Tc\|^2}{\|c\|_{1,w}^2}\end{equation}
 Indeed,\begin{align*}
 \|Tc\|^2&=\langle Tc, \sum_{j=1}^N c_jf_j\rangle\\
         &=\sum_{j=1}^{N}c_j\langle Tc,f_j\rangle\\
         &\leq\sum_{j=1}^{N}|c_j||w_j|\frac{|\langle Tc,f_j\rangle|}{\|f_j\|}\\
         &\leq\frac{|\langle Tc,f_{j_0}\rangle|}{\|f_{j_0}\|}\|c\|_{1,w}.
 \end{align*}
But $j_0\notin\supp c$ implies that $|\langle f_j,f_{j_0}\rangle|\leq\mu\|f_j\|\|f_{j_0}\|,$ for $j\in\supp c.$\\
 Hence $$|\langle Tc,f_{j_0}\rangle|=|\sum_{j\in\supp c}c_j\langle f_j,f_{j_0}\rangle|\leq\mu\sum_{j\in\supp c}|c_j|\|w_j\|\|f_{j_0}\|=\mu\|f_{j_0}\|\|c\|_{1,w}.$$
 Hence $$\|Tc\|^2\leq\mu\|c\|_{1,w}^2.$$
 Using Cauchy-Schwarz inequality, from relation (\ref{suportul}), we have
 \begin{align*}
 \mu&\geq\frac{\|Tc\|^2}{\|c\|_{1,w}^2}\geq\frac{(1-\delta_s)\|c\|_{2,w}^2}{\|c\|_{1,w}^2}\\
    &\geq\frac{1-\delta_s}{s}.
  \end{align*}
 \end{pf}
 \begin{cor}
 If $c$ is $s-$sparse and $\mu<\dfrac{1}{2s-1},$ then $j_0\in\supp c.$
 \end{cor}
 \begin{pf}
 From Proposition \ref{prop1sup} and Theorem \ref{thm2sup}.
 \end{pf}
\section{Main results}
For a vector $c\in\mathbb{R}^N$ we denote by $c_s$ the vector $c$ with all but the $s-$largest entries set to zero. Also, for $T_0\subset\{1,2,\ldots,N\}$, we denote by $T_0^c$ the complement of $T_0$.
\begin{thm}\label{Thm5}
Let $c,d\in\mathbb{R}^N$ such that $\|c\|_{1,w}\geq\|d\|_{1,w}$ and let $v=d-c$. If $$\mu<\frac{1}{2s-1}$$ then we have $$\|v\|_{2,w}\leq\frac{\sqrt{3-\frac{1}{2s-1}}}{1-\mu(2s-1)}\|Tv\|+\frac{2\sqrt{\mu(1+\mu)s}}{1-\mu(2s-1)}e_0,$$
where $$e_0=\frac{\|c-c_s\|_{1,w}}{\sqrt{s}}.$$
\end{thm}
\begin{pf}
We denote by $T_0$ the locations of the $s$ largest coefficients of $c$. We have
\begin{equation}\label{eq0}
\|v_{T_0^c}\|_{1,w}\leq\|v_{T_0}\|_{1,w}+2e_0.
\end{equation}
Indeed,
\begin{align*}
\|c\|_{1,w}&\geq\|d\|_{1,w}=\sum_{i=1}^Nw_i|v_i+c_i|\\
           &=\sum_{i\in T_0}w_i|v_i+c_i|+\sum_{i\in T_0^c}w_i|v_i+c_i|\\
           &\geq\|c_{T_0}\|_{1,w}-\|v_{T_0}\|_{1,w}+\|v_{T_0^c}\|_{1,w}-\|c_{T_0}\|_{1,w},
\end{align*}
hence $$2\|c_{T_0^c}\|_{1,w}+\|v_{T_0}\|_{1,w}\geq\|v_{T_0^c}\|_{1,w},$$ i.e. relation (\ref{eq0}). We have
\begin{align*}
\langle Tv,Tv_{T_0}\rangle&=\langle Tv_{T_0}+Tv_{T_0^c},Tv_{T_0}\rangle\\
                          &=\|Tv_{T_0}\|^2+\langle Tv_{T_0^c},Tv_{T_0}\rangle\\
                          &\geq[1-\mu(s-1)]\|v_{T_0}\|_{2,m}^2-\mu\|v_{T_0^c}\|_{1,m}\|v_{T_0}\|_{1,w}
\end{align*}
But $v_{T_0}$ is $s-$ sparse and from Lemmma \ref{Basic lemma} we have
\begin{equation}\label{eq1}
\|v_{T_0}\|_{1,w}\leq\sqrt{s}\|v_{T_0}\|_{2,w}
\end{equation}
We obtain $$\|Tv\|\|Tv_{T_0}\|\geq[1-\mu(s-1)]\|v_{T_0}\|_{2,w}^2-\mu\|v_{T_0^c}\|_{1,m}\sqrt{s}\|v_{T_0}\|_{2,w}$$
But
\begin{equation}\label{eq2}
\|Tv_{T_0}\|\leq\sqrt{1+\mu(s-1)}\|v_{T_0}\|_{2,w}
\end{equation}
So $$\sqrt{1+\mu(s-1)}\|Tv\|\geq[1-\mu(s-1)]\|v_{T_0}\|_{2,w}-\mu\sqrt{s}\|v_{T_0^c}\|_{1,w}$$
But
\begin{equation}\label{eq3'}
\|v_{T_0^c}\|_{1,w}\leq\|v_{T_0}\|_{1,w}+2\sqrt{s}e_0.
\end{equation}
Using relations (\ref{eq1}) and (\ref{eq3'}), we obtain
\begin{equation}\label{eq3}
\|v_{T_0^c}\|_{1,w}\leq\sqrt{s}[\|v_{T_0}\|_{2,w}+2e_0]
\end{equation}
It follows that
\begin{align*}\sqrt{1+\mu(s-1)}\|Tv\|&\geq[1-\mu(s-1)]\|v_{T_0}\|_{2,w}-\mu\sqrt{s}[\sqrt{s}(\|v_{T_0}\|_{2,w}+2e_0)]\\
                                     &\geq[1-\mu(s-1)]\|v_{T_0}\|_{2,w}-\mu s\|v_{T_0}\|_{2,w}-2\mu se_0
\end{align*}
So we obtain the following inequality
\begin{equation}\label{eq4}
\|v_{T_0}\|_{2,w}\leq\frac{\sqrt{1+\mu(s-1)}}{1-\mu(2s-1)}\|Tv\|+\frac{2\mu se_0}{1-\mu(2s-1)}
\end{equation}
On the other hand, by the basic lemma (Lemma \ref{Basic lemma}), we have
$$\|Tv\|^2\geq(1+\mu)\|v\|_{2,w}^2-\mu\|v\|_{1,w}^2$$
and by using the equation (\ref{eq1}) and the equation (\ref{eq3'}) we obtain
\begin{align*}
\|v\|_{1,w}&=\|v_{T_0}\|_{1,w}+\|v_{T_0^c}\|_{1,w}\\
           &\leq 2\|v_{T_0}\|_{1,w}+2\sqrt{s}e_0\\
           & \leq 2\sqrt{s}\|v_{T_0}\|_{2,w}+2\sqrt{s}e_0.
\end{align*}
so \begin{equation}\label{eq5}
\|v\|_{1,w}\leq 2\sqrt{s}\|v_{T_0}\|_{2,w}+2\sqrt{s}e_0.
\end{equation}
and then $$\|Tv\|^2\geq(1+\mu)\|v\|_{2,w}^2-4\mu s(\|v_{T_0}\|_{2,w}+e_0)^2$$
which is equivalent to
\begin{equation}\label{eq6}
(1+\mu)\|v\|_{2,w}^2\leq\|Tv\|^2+4\mu s(\|v_{T_0}\|_{2,w}+e_0)^2
\end{equation}
Combining the equation (\ref{eq4}) and the equation (\ref{eq6}), we have
\begin{align*}
(1+\mu)\|v\|_{2,w}^2&\leq\|Tv\|^2+4\mu s\bigg(\frac{\sqrt{1+\mu(s-1)}}{1-\mu(2s-1)}\|Tv\|+\frac{2\mu se_0}{1-\mu(2s-1)}+e_0\bigg)^2\\
&=\|Tv\|^2+4\mu s\bigg(\frac{\sqrt{1+\mu(s-1)}}{1-\mu(2s-1)}\|Tv\|+\frac{(1+\mu) e_0}{1-\mu(2s-1)}\bigg)^2.
\end{align*}
Hence
$$ (1+\mu)\|v\|_{2,w}^2 \leq \|Tv\|^2+\bigg(2\sqrt{\mu s}\frac{\sqrt{1+\mu(s-1)}}{1-\mu(2s-1)}\|Tv\|+2\sqrt{\mu s}\frac{(1+\mu )e_0}{1-\mu(2s-1)}\bigg)^2.$$

We use the following inequality $$\alpha^2+(m\alpha+\beta)^2\leq(\sqrt{1+m^2}\alpha+\beta)^2$$
where
\begin{align*}\alpha&=\|Tv\|\\
m&=2\sqrt{\mu s}\frac{\sqrt{1+\mu(s-1)}}{1-\mu(2s-1)}\\
\beta&=2\sqrt{\mu s}\frac{(1+\mu) e_0}{1-\mu(2s-1)}
\end{align*}
With this notations we have $$(\sqrt{1+\mu})^2\|v\|_{2,w}^2\leq(\sqrt{1+m^2}\alpha+\beta)^2$$
so $$\|v\|_{2,w}^2\leq\frac{1}{\sqrt{1+\mu}}(\sqrt{1+m^2}\alpha+\beta)^2$$
But \begin{align*}
\sqrt{1+m^2}=\frac{\sqrt{(8s^2-8s+1)\mu^2+2\mu+1}}{1-\mu(2s-1)}
\end{align*}
So we obtain $$\|v\|_{2,w}\leq\frac{1}{\sqrt{1+\mu}}\bigg( \frac{\sqrt{(8s^2-8s+1)\mu^2+2\mu+1}}{1-\mu(2s-1)}\|Tv\|+\frac{2\sqrt{\mu s}(1+\mu )e_0}{1-\mu(2s-1)}\bigg)$$
And so $$\|v\|_{2,w}\leq\frac{\|Tv\|}{1-\mu(2s-1)}\sqrt{\frac{(8s^2-8s+1)\mu^2+2\mu+1}{1+\mu}}+\frac{2\sqrt{\mu s}\sqrt{1+\mu }e_0}{1-\mu(2s-1)}.$$
We denote $$F(\mu)=\frac{(8s^2-8s+1)\mu^2+2\mu+1}{1+\mu}$$
Since $$F'(\mu)=\frac{(8s^2-8s+1)\mu^2+2(8s^2-8s+1)\mu+1}{(1+\mu)^2}\geq 0.$$ and $\mu\leq\frac{1}{2s-1}$ it follows that $$F(\mu)\leq F(\frac{1}{2s-1})=\frac{6s-4}{2s-1}.$$
Finally, we obtain that $$\|v\|_{2,w}\leq\frac{\sqrt{3-\frac{1}{2s-1}}}{1-\mu(2s-1)}\|Tv\|+\frac{2\sqrt{\mu s(1+\mu)}}{1-\mu(2s-1)}e_0.$$
\end{pf}
\begin{thm} Assume that $$\mu<\frac{1}{2s-1}$$ and $\|z\|_2\leq\varepsilon$. Then the solution $c^*$ of $(P_{1,w})$ obeys
$$\|c^*-c\|_{2,w}\leq\frac{\sqrt{3-\frac{1}{2s-1}}}{1-\mu(2s-1)}(\eta+\varepsilon)+\frac{2\sqrt{\mu s(1+\mu)}}{1-\mu(2s-1)}e_0.$$
\end{thm}
\begin{pf}
This result it follows from Theorem \ref{Thm5} since
\begin{align*}
\|Tv\|_2&=\|Tc^*-y-(Tc-y)\|_2\\
        &\leq\|Tc^*-y\|_2+\|Tc-y\|_2\\
        &\leq\varepsilon+\eta.
\end{align*}
\end{pf}
In \cite{Cai} the authors obtained, for $c$ $s-$sparse and $\|f_i\|=1, i=\overline{1,N}$, the following estimation
$$\|c^*-c\|_{2,w}\leq\frac{\sqrt{3(1+\mu)}}{1-(2s-1)\mu}(\eta+\varepsilon).$$ We notice that our estimation, given in the above Theorem is better. The authors of the paper presented in reference \cite{Cai} indicate that there exist a relation like the one given in the above Theorem, but without the specification of the constants.

\end{document}